\documentclass[10pt]{article}

\usepackage[cp1250]{inputenc}
\usepackage[english]{babel}

\usepackage{amsfonts,amsmath,amssymb,amsthm}
\usepackage[abbrev]{amsrefs}

\usepackage{graphicx}
    \DeclareGraphicsExtensions{.pdf}

\usepackage{tikz}
\usetikzlibrary{matrix}

\newtheoremstyle{ex}
{3pt}
{3pt}
{}
{}
{}
{.}
{.5em}
{}

\theoremstyle{definition}
\newtheorem{thm}{Theorem}[section]
\newtheorem{prop}[thm]{Proposition}
\newtheorem{cor}[thm]{Corollary}
\newtheorem{lemma}[thm]{Lemma}

\theoremstyle{ex}
\newtheorem{dfn}[thm]{Definition}
\newtheorem{example}[thm]{Example}
\newtheorem{rem}[thm]{Remark}
\newtheorem*{pf}{Proof}

\newcommand{\z}[1]{\mathbb{Z}/#1}

\author{Wojciech Politarczyk}
\title{$4$-manifolds, surgery on loops and geometric realization of Tietze transformations}

\begin{document}
\maketitle

\begin{abstract}
In the paper \cite{wall_1}, C.T.C. Wall proved that two smooth closed simply connected $4$-manifolds which are homeomorphic are in fact stably diffeomorphic. We prove a similar result which states that two smooth closed $4$-manifolds satisfying certain properties are stably diffeomorphic if and only if their signatures agree. The manifolds in question are obtained by surgery on loops. The methods we use are modified surgery of Kreck \cite{kreck} and Kirby calculus.
\end{abstract}

\let\thefootnote\relax

\section{Introduction}
A well known theorem states that for any finitely presented group $\pi$ and $n \geq 4$ there exists a smooth compact connected $n$-manifold whose fundamental group is isomorphic to $\pi$. The construction starts with arbitrary presentation $\mathcal{P}$ of $\pi$. Then the surgery is performed on
$$(S^1 \times S^{n-1}) \# \ldots \# (S^1 \times S^{n-1})$$
along loops representing the relations of $\mathcal{P}$. The Seifert-Van Kampen theorem implies that the obtained manifold has fundamental group isomorphic to $\pi$.

The manifold obtained by the construction will depend on the presentation of $\pi$ and the way the surgery was performed. In order to classify these manifolds one has to take care of all the data, improve the construction and make it more precise. This is done in section \ref{section_surgery_on_loops} for $n=4$. In this section the construction of $\pi$-realizing surgery on loops is defined.

In 1960's C.T.C. Wall proved in \cite{wall_1} that two smooth closed simply-connected $4$-manifolds which are h-cobordant become diffeomerphic after taking the connected sum with a number of copies of the product $S^2 \times S^2$. This leads to the notion of stable diffeomorphism.
\begin{dfn}
Let $M$ and $N$ be two smooth manifolds of dimension $2n$. We say $M$ and $N$ are \textit{stably diffeomorphic} if there are two positive integers $k,l$ and a diffeomorphism
$$f \colon M \# k(S^n \times S^n) \to N \#l (S^n \times S^n).$$ 
\end{dfn}
The stable diffeomorphism relation played an essential role in the later developments of the $4$-manifold topology. For example from the results of Wall and Milnor, contained in \cite{wall_1} and \cite{milnor}, one can deduce that any two closed smooth simply connected $4$-manifolds with the same homotopy type are stably diffeomorphic. In particular any closed smooth simply connected $4$-manifold has a unique smooth structure up to stable diffeomorphism. This results was generalized by other authors. For example Gompf in \cite{gompf} proved that the conclusion about the stable uniquness of a smooth structure holds for any fundamental group. J.F. Davis proved in \cite{davis} that if a group $\pi$ satisfies certain algebraic condition then homotopy equivalence implies stable diffeomorphism for manifolds with fundamental group isomorphic to $\pi$. It is also worth mentioning that Quinn proved in \cite{quinn} that s-cobordant $4$-manifolds are stably diffeomorphic.


The main result of the first part of this paper reads as follows (for the definition of normal $1$-type see def. \ref{dfn_normal_1_type}):

\begin{thm}\label{main_theorem}
Let $\pi$ be a finitely presented group. Two closed connected orientable smooth $4$-manifolds obtained by $\pi$-realizng surgery on loops with the same normal $1$-type are stably diffeomorphic if and only if the signatures of both manifolds are equal.
\end{thm}


\begin{example}
Suppose that we chose two finite presentations $\mathcal{P}_1$ and $\mathcal{P}_2$ of a fixed group $\pi$. Let $M_1 = \#{k_1}(S^3 \times S^1)$ and $M_2 = \#{k_2} (S^3 \times S^1)$, where $k_1$ and $k_2$ are the numbers of generators in $\mathcal{P}_1$ and $\mathcal{P}_2$, respectively. 

The normal $1$-type of even $4$-manifold with fundamental group $\pi$ is determined by a cohomology class $w \in H^2(\pi ; \z{2})$. Two classes $w_1, w_2 \in H^2(\pi ; \z{2})$ define the same normal $1$-type if there exists an automorphism $\phi \colon \pi \to \pi$ such that
$$\phi^{\ast}w_1 = w_2.$$
For more details see the discussion after the lemma \ref{lemma_exact_seq}.

In our case for any cohomology class $w \in H^2(\pi ; \z{2})$ the $\pi$-realizing surgery on loops can be performed on $M_1$ and $M_2$ such that the normal $1$-type corresponds to~$w$.

Let $w_1,w_2 \in H^2(\pi ; \z{2})$ and $Y(M_1,\mathcal{P}_1,w_1)$, $Y(M_2,\mathcal{P}_2,w_2)$ denote the result of the $\pi$-realizing surgery on loops performed on $M_1$ and $M_2$, respectively. The normal $1$-type of $M_i$ is corresponds to the class $w_i$ for $i=1,2$. Theorem \ref{main_theorem} implies, due to the equality of signatures, that $Y(M_1,\mathcal{P}_2,w_1)$ and $Y(M_2,\mathcal{P}_2,w_2)$ are stably diffeomorphic if and only if there exists an automorphism $\phi$ of $\pi$ such that $$\phi^{\ast}w_1 = w_2.$$
In particular if $w_1=w_2=0$, then both manifolds are stably diffeomorphic.
\end{example}

\begin{example}
Let $M_1' = M_1 \# \mathbb{C}P^2$ and $M_2' = M_2 \# \mathbb{C}P^2$, where $M_1$ and $M_2$ are the same as in the previous example. For odd $4$-manifolds the normal $1$-type is determined by the fundamental group. Thus if $Y(M_1',\mathcal{P}_1)$ and $Y(M_2',\mathcal{P}_2)$ are always stably diffeomorphic.
\end{example}

In the second part of the paper the considerations will focus on manifolds obtained from the connected sums of $S^1 \times S^3$. The aim of this part is to reprove the special case of the theorem \ref{main_theorem}in two steps. The first step consists of showing that any  manifold obtained from the connected sum of copies of $S^3 \times S^1$ by $\pi$-realizing surgery on loops, is diffeomorphic to a thickening $T(X,\xi)$ of a $2$-dimensional CW complex $X$. For more informations about thickenings see section \ref{section_thickenings}.

\begin{thm}\label{thm_diffeo_thick_surg}
Let $\mathcal{P}$ be a presentation of a finitely presented group $\pi$ and let $w \in H^2(\pi,\z{2})$ be arbitrary. Let $Y(X_n,\mathcal{P},w)$ denote a result of a $\pi$-realizing surgery on loops performed on a connected sum $\#_n S^1 \times S^3$. Then there exists a finite $2$-dimensional CW complex $X$ and a vector bundle $\xi$ over it such that $Y(X_n,\mathcal{P},w)$ is diffeomorphic to $T(X,\xi)$.
\end{thm}

The second step uses the fact that for thickenings an appropriate version of theorem \ref{main_theorem} can be proved using Kirby calculus. The starting point is the fact that every $2$-dimensional CW-complex correspond to a presentation of its fundamental group. The next step is to use Tietze theorem. This theorem states that two presentations yield the same groups if and only if one of them can be transformed into the other by a finite sequence of Tietze transformations. In section \ref{section_geometric_realization} geometric realizations of Tietze transformation are described. These geometric transformation are sequences of handle slides performed on Kirby diagrams of thickenings $T(X,\xi)$. If $\mathcal{P}_1$ and $\mathcal{P}_2$ are two presentations of the same group, differing by an application of a single Tietze transformations, then performing appropriate geometric transformation on $T(X(\mathcal{P}_2),\xi_1)$ yields $T(X(\mathcal{P}_2),\xi_2)$. Thorough analysis of geometric Tietze transformations 
yields the following theorem.

\begin{thm}\label{thm_geometric_realisation_stable_diffeo}
Let $\mathcal{P}_1$ and $\mathcal{P}_2$ be two presentations of a fixed group $\pi$ differing by a single application of a Tietze transformation. Let $X_1$ and $X_2$ be two finite $2$-dimensional CW-complexes corresponding to $\mathcal{P}_1$ and $\mathcal{P}_2$. Let $\xi_1$ and $\xi_2$ be vector bundles over $X_1$ and $X_2$ of dimension $3$. If the thickenings have the same normal $1$-type, then two cases can occur. Either there exists a diffeomorphism
$$T(X_1,\xi_1) \# S^2 \times S^2 \stackrel{\cong}{\longrightarrow} T(X_2,\xi_2)$$
or thickenings $T(X_1,\xi_1)$ and $T(X_2,\xi_2)$ are diffeomorphic.
\end{thm}

The next theorem is an easy corollary of the theorem \ref{thm_geometric_realisation_stable_diffeo}.

\begin{thm}\label{thm_main_1}
Let $X_1$ and $X_2$ be two finite $2$-dimensional CW-complexes with isomorphic fundamental groups. Let $\xi_1$ and $\xi_2$ be two  vector bundles of dimension three over $X_1$ and $X_2$, respectively. Then the thickenings $T(X_1,\xi_1)$ and $T(X_2,\xi_2)$ are stably diffeomorphic if they have the same normal $1$-type.
\end{thm}

The paper is divided into two parts. The first part consists of section \ref{section_normal_types}, \ref{section_surgery_on_loops} and \ref{section_stable_classification}. The second part consists of sections \ref{section_thickenings}, \ref{section_Tietze}, \ref{section_geometric_realization} and \ref{section_surgery_kirby_diagram}.

In section \ref{section_normal_types} the normal $1$-type is defined. This concept is essential for controlling changes made during the construction. It is also a starting point of the modified surgery of Kreck, see \cite{kreck}. The advantage of using the normal $1$-type comes from the Kreck's stable diffeomorphism theorem (thm. \ref{krecks_thm}). One reduces the problem of stable diffeomorphism classification to the algebraic problem of computing the bordism groups of the fibration representing the normal $1$-type.

Section \ref{section_surgery_on_loops} describes the surgery on loops under circumstances needed in this paper.

Section \ref{section_stable_classification} contains the proof of the main theorem \ref{main_theorem}.

Section \ref{section_thickenings} is concerned with the construction and properties of thickenings of CW-complexes. In this section the Kirby diagrams for such thickenings are constructed.

In section \ref{section_Tietze} a brief recollection of Tietze transformations and Tietze theorem is provided.

Secion \ref{section_geometric_realization} is concerned with the description of geometric Tietze transformations. These transformations are used to prove theorem \ref{thm_geometric_realisation_stable_diffeo}.

Section \ref{section_surgery_kirby_diagram} contains description of changes that a Kirby diagram undergoes during the surgery. This description is later used to prove theorem \ref{thm_diffeo_thick_surg}.

From now on every manifold under consideration is smooth compact and connected.

\section{Normal 1-types}\label{section_normal_types}
Presentasion of the material in this section follows \cite{kreck}.

\begin{dfn}\label{dfn_normal_1_type}
Let $k \in \mathbb{N}$ and  $p \colon B_{k} \to BSO$ be a fibration. We say that $p$ is \textit{$k$-universal} if the homotopy fibre of $p$ is connected and the relative homotopy groups $\pi_i(p)$ vanish for $i \geq k+1$.

\noindent Let $M$ be a manifold with a $B_k$-structure, i.e. the normal Gauss map $$\nu \colon M \to BSO$$ admits a lifting $\overline{\nu}$ to $B_k$. The map $\overline{\nu}$ is called \textit{$k$-smoothing} if it is a $(k+1)$-equivalence.
\end{dfn}

\begin{dfn}
Let $M$ be a closed connected and oriented $n$-manifold. The normal \textit{$k$-type of $M$} is the fibre homotopy type of the $k$-universal fibration
$$B_k \to BSO,$$
such that there exists a $k$-smoothing of $M$
\begin{center}
\begin{tikzpicture}
\matrix(m)[matrix of math nodes, column sep=1cm,row sep=1cm]
{
 & B_k \\
M & BSO \\
};
\path[->,font=\scriptsize]
(m-2-1) edge node[above] {$\nu$} (m-2-2)
(m-1-2) edge node {} (m-2-2);
\path[->,dashed,font=\scriptsize]
(m-2-1) edge node[left] {$\overline{\nu}$} (m-1-2);
\end{tikzpicture}
\end{center}
\end{dfn}
In this paper we will be mostly interested in normal $1$-types of $4$-manifolds.
To construct a fibration, which represents the normal $1$-type of a given $4$-manifold~$M$ we have to consider two separate cases.

First suppose that $w_2(\widetilde{M}) \neq 0$, where $\widetilde{M}$ is the universal covering space of~$M$ and $w_2$ is the second Stiefel-Whitney class.
\begin{dfn}\label{dfn_odd_type}
A $4$-manifold $M$ is said to be \textit{odd} if its universal covering space is not a spin manifold.
\end{dfn}
In this case the normal $1$-type of $M$ is represented by the following fibration.
\begin{center}
\begin{tikzpicture}
\matrix(m)[matrix of math nodes, column sep=1cm,row sep=1cm]
{
 & BSO \times K(\pi_1(M),1)\\
M & BSO \\
};
\path[->,font=\scriptsize]
(m-2-1) edge node[above] {$(\nu,u)$} (m-1-2)
		  edge node[above] {$\nu$} (m-2-2)
(m-1-2) edge node {} (m-2-2);
\end{tikzpicture}
\end{center} 
Here $u \colon M \to K(\pi_1(M),1)$ is the classifying map for the universal covering.

\begin{dfn}\label{dfn_even_type}
A $4$-manifold is said to be \textit{even} if its universal covering space is a spin manifold.
\end{dfn}
Let $M$ be even and let $u \colon M \to K(\pi_1(M),1)$ classify the universal covering. Consider the following lemma.
\begin{lemma}\label{lemma_exact_seq}
Suppose $M$ is a finite CW-complex and let $A$ be a trivial $\mathbb{Z}\pi_1(M)$ module. Then the following sequence is exact
$$0 \to H^2(\pi_1(M),A) \stackrel{u^{\ast}}{\to} H^2(M,A) \stackrel{p^{\ast}}{\to} H^2(\widetilde{M},A).$$
\end{lemma}

From the lemma \ref{lemma_exact_seq} we can conclude that there exists a unique class $w \in H^2(\pi,\z{2})$ such that $u^{\ast}w=w_2(M)$. Consider now the following pullback square.
\begin{center}
\begin{tikzpicture}
\matrix(m)[matrix of math nodes, row sep=1cm, column sep=1cm]
{
B(\pi,w) & K(\pi,1) \\
BSO & K(\z{2},2)\\
};
\path[->,font=\scriptsize]
(m-1-1) edge node {} (m-1-2)
(m-1-1) edge node {} (m-2-1)
(m-1-2) edge node[right] {$w$} (m-2-2)
(m-2-1) edge node[above] {$w_2$} (m-2-2);
\end{tikzpicture}
\end{center}
\begin{prop}\label{prop_even_normal_1_type}
The fibration $B(\pi,w) \to BSO$ from the diagram above represents the normal $1$-type of $M$.
\end{prop}
\begin{pf}
Two maps $u$ and $\nu$ determine the lifting $\overline{\nu} \colon M \to B(\pi,w)$ of the normal map $\nu$ by the universal property of the fibration. Observe now that $BSpin$ is the homotopy fibre of the map $B(\pi,w) \to K(\pi,1)$, which implies that this map induces an isomorphism on $\pi_1$ and $\pi_2$. The map $u$ is a $2$-equivalence, thus $\overline{\nu}$ must be a $2$-equivalence. Additionally the homotopy groups of the homotopy fibre of the map $B(\pi,w) \to BSO$ vanish in dimension greater than $1$. This proves that $B(\pi,w) \to BSO$ is a $1$-universal fibration and $\overline{\nu}$ is a $2$-smoothing.
\end{pf}

\begin{rem}
The proposition \ref{prop_even_normal_1_type} implies that for even $4$-manifolds the normal $1$-type is determined by the fundamental group $\pi$  and a certain cohomology class $w \in H^2(\pi;\z{2})$. Is is easy to check that two pairs $(\pi_1,w_1)$ and $(\pi_2,w_2)$ determine fibre homotopy equivalent fibrations if and only if there exists an isomorphism $\phi \colon \pi_1 \to \pi_2$ such that $\phi^{\ast}w_2 = w_1$.
\end{rem}

\section{$\pi$-realizing surgery on loops}\label{section_surgery_on_loops}
First suppose that we are given a pair $(\pi_1, w_1)$, where $\pi_1$ is a finitely presented group and $w_1 \in H^2(\pi_1,\z{2})$. Using this data one can construct a $1$-universal fibration over $BSO$ such that the following diagram is a pullback diagram.
\begin{center}
\begin{tikzpicture}
\matrix(m)[matrix of math nodes, row sep=1cm, column sep=1cm]
{
B(\pi_1,w_1) & K(\pi_1,1) \\
BSO & K(\z{2},2)\\
};
\path[->,font=\scriptsize]
(m-1-1) edge node {} (m-1-2)
(m-1-1) edge node {} (m-2-1)
(m-1-2) edge node[right] {$w$} (m-2-2)
(m-2-1) edge node[above] {$w_2(\gamma)$} (m-2-2);
\end{tikzpicture}
\end{center}

\noindent Let $(\pi_2,w_2)$ be another such pair and suppose that there exists an epimorphism $\psi \colon \pi_1 \to \pi_2$ such that $w_1 = \psi^{\ast} w_2$. Then by the universal property of the pullback there exists a map of fibrations
$$B(\pi_1,w_1) \stackrel{\psi_{\#}}{\to} B(\pi_2,w_2)$$
Suppose now that $M$ is a smooth, closed, oriented, even $4$-manifold whose normal $1$-type is fibre-homotopy equivalent to $B(\pi_1,w_1) \to BSO$. Let $$\overline{\nu} \colon M \to B(\pi_1,w_1)$$ be the lift of the normal Gauss map. Then the composition $$\psi_{\#} \circ \overline{\nu} \colon M \to B(\pi_2,w_2)$$ is a normal map. We can perform surgery below the middle dimension on this map to a obtain $1$-smoothing $Y(M,\psi,w_2) \to B(\pi_2,w_2)$.

\begin{dfn}
We say that the manifold $Y(M,\psi,w_2)$ was obtained from $M$ by \textit{$\pi_2$-realizing surgery on loops of type I}.
\end{dfn}

\begin{prop}\label{prop_homomorphism_spin}
Surgery on loops of type I yields a homomorphism of bordism groups $\psi_{\ast} \colon \Omega_{4}(B(\pi_1,w_1)) \to \Omega_4(B(\pi_2,w_2))$
\end{prop}

\begin{rem}
To be precise one has to take care of the framings to be sure that the obtain manifold has the desired normal $1$-type. If $\gamma \colon S^1 \to M$ is an embedding whose homotopy class represents one of the relations, then the framing of the stable normal bundle of the loop $\gamma$ comes from the nullhomotopy $H$ of the composition $\psi_{\#} \circ \overline{\nu} \circ \gamma$. If we take different nullhomotopy $H'$, then these two maps determine a map $S^2 \to B(\pi_1,w_1)$. However it is easy to check that $\pi_2(B(\pi_1,w_1))=0$, which implies that any two nullhomotopies are homotopic relative to the boundary, thus the homotopy class of the stable framing is uniquely defined.
\end{rem}

In the odd case the surgery construction is analogous. Let $Y(M,\psi)$ be the result of the construction performed on $M$ using the homomorphism $\psi$.

\begin{dfn}
We say that $Y(M,\psi)$ was obtained from $M$ by \textit{$\pi_2$-realizing surgery on loops of type II}.
\end{dfn}

\begin{prop}\label{prop_homomorphism_nonspin}
If $\psi \colon \pi_1 \to \pi_2$ is an epimorphism, then the surgery on loops yields a homomorphism of bordism groups
$$\psi_{\ast} \colon \Omega_4^{SO}(K(\pi_1,1)) \to \Omega_4^{SO}(K(\pi_2,1)).$$
\end{prop}

\section{Stable classification of $4$-manifolds} \label{section_stable_classification}
\noindent Using Kreck's modified surgery from \cite{kreck} one can give a sufficient condition for two manifolds $M$ and $N$ to be stably diffeomorphic.
\begin{thm}[\cite{kreck}]\label{krecks_thm}
Suppose $M$ and $N$ are two smooth, connected and closed $2q$-manifolds with the same normal $(q-1)$-type $B^{q-1} \to BO$, $q \geq 2$. If $(M, \overline{\nu}_{M}) = (N, \overline{\nu}_{N})$ in $\Omega_{2q}^{B^{q-1}}$, then $M$ and $N$ are stably diffeomorphic.
\end{thm}
\noindent Theorem \ref{krecks_thm} yields the following corollary.
\begin{cor}
Let $M_1$ and $M_2$ be two $4$-manifolds whose normal $1$-type is fibre-homotopy equivalent to $B(\pi_1,w_1) \to BSO$. Let $\psi \colon \pi_1 \to \pi_2$ be an epimorphism of finitely presented groups and $w_2 \in H^2(\pi_2,\z{2})$ such that $w_1 = \psi^{\ast} w_2$. Then $Y(M_1,\psi,w_2)$ and $Y(M_2,\psi,w_2)$ are stably diffeomorphic if and only if $$[M_1, \overline{\nu_{M_1}}] - [M_2, \overline{\nu_{M_2}}] \in \ker(\psi_{\ast} \colon \Omega_4(B(\pi_1,w_1)) \to \Omega_4(B(\pi_2,w_2))).$$
\end{cor}
\noindent In order to prove the main theorem one has to investigate the properties of the homomorphisms from propositions \ref{prop_homomorphism_spin} and \ref{prop_homomorphism_nonspin}.

The proof of theorem \ref{main_theorem} requires some knowledge of the bordism group $\Omega_4(B,\pi,w)$ of the $1$-universal fibration $B(\pi,w) \to BSO$. When $w=0$ this group is isomorphic to the spin bordism group of the Eilenberg-MacLane space $K(\pi,1)$ thus the classical Atiyah-Hirzebruch spectral sequence can be used to extract the necessary information. To deal with the case when $w \neq 0$ we need the James spectral sequence, which was constructed by Teichner in \cite{teichner}.

\begin{thm}[James spectral sequence (JSS)]\label{thm_jss}
Let $F \to B \to K$ be a fibration and $B \to BSO$ be a stable vector bundle. Suppose that $h$ is a generalized homology theory whose coefficients vanish in negative dimensions. Then there exists a spectral sequence
$$E^2_{p,q} = H_{p}(K,h_q(Th(\xi|_{F}))) \Longrightarrow h_{p+q}(Th(\xi)),$$
where $Th$ stands for the Thom space of apropriate bundle.

Furthermore this spectral sequence is natural with respect to the maps of fibrations.
\end{thm}

\begin{prop}
Let $\pi$ be a finitely presented group. Let $\psi \colon F_{\alpha} \to \pi$ be an epimorphism. Then the homomorphism $\psi_{\ast} \colon \Omega_4^{Spin}(K(F_{\alpha},1)) \to \Omega_4(B(\pi,w))$ is injective. For any two epimorphisms $\psi_i \colon F_{\alpha_i} \to \pi$, $i=1,2$, we have $$\operatorname{im}\left((\psi_1)_{\ast} \right) = \operatorname{im}\left( (\psi_2)_{\ast} \right).$$
\end{prop}
\begin{pf}
The proof uses naturality of JSS. First observe that $H_p(F_{\alpha})=0$ for $p\geq 2$ and $\Omega_3^{Spin}=0$. Thus $E^2_{1,3}=E^2_{2,2}=E^2_{3,1}=E^2_{4,0}=0$. Naturality implies that $\operatorname{im}(\psi_{\ast}) \subset E^{\infty}_{0,4} = E^2_{0,4} = \Omega_4^{Spin}$. The edge homomorphism corresponding to this entry is given by the signature.
\end{pf}

\noindent As a corollary we obtain the following theorem.
\begin{thm}\label{classification_even_type}
Let $\pi$ be a finitely presented group and $w \in H^2(\pi,\z{2})$. Choose two presentations of $\pi$ and let $\psi_1 \colon F_{k_1} \to \pi$ and $\psi_2 \colon F_{k_2} \to \pi$ be homomorphisms induced by these presentations. Choose two spin $4$-manifolds $M_1$ and $M_2$ such that $\pi_1(M_1) = F_{k_1}$ and $\pi_1(M_2)=F_{k_2}$. Then $Y(M_1,\psi_{1},w)$ and $Y(M_2,\psi_{2},w)$ are stably diffeomorphic if and only if $\sigma(M_1)=\sigma(M_2)$.
\end{thm}

\begin{rem}
It is natural to ask whether every even $4$-manifold with fundamental group $\pi$, is stably diffeomorphic to a one obtained from a spin manifold with free fundamental group. The answer to this question is negative. Rochlin's theorem implies that necessarily the signature should be divisible by $16$. According to the main theorem from \cite{teichner2} if $H^2(\pi;\z{2}) \neq 0$ then there exists an even $4$-manifold with fundamental group isomorphic to $\pi$ and signature equal to~$8$. Even when $H^2(\pi;\z{2})=0$ there are other obstructions coming from the other bordism invariants.
\end{rem}

\noindent In the odd case we have the following lemma.
\begin{lemma}
The following map is an isomorphism
\begin{align*}
\Omega_4^{SO}(K(\pi,1)) &\to \mathbb{Z} \oplus H_4(K(\pi,1);\mathbb{Z})\\
(M,f) &\mapsto (\sigma(M),f_{\ast}[M]).
\end{align*}
\end{lemma}
\begin{pf}
The lemma follows easily from the investigation of the Atiyah-Hirzebruch spectral sequence. Indeed notice that $\Omega_1^{SO} = \Omega_2^{SO}=\Omega_3^{SO}=0$ implies that $E^2_{1,3} = E^2_{2,2} = E^2_{3,1}=0$. Thus there are no nonzero differentials. Therefore we obtain the following exact sequence.
$$0 \to \Omega_4^{SO} \to \Omega_4^{SO}(K(\pi,1)) \to H_4(\pi ; \mathbb{Z}) \to 0.$$
However the constant map $K(\pi,1) \to \ast$ induces a map on the bordism groups which splits the inclusion $\Omega_4^{SO} \to \Omega_4^{SO}(K(\pi,1))$. Also one can show that the epimorphism $\Omega_4^{SO}(K(\pi,1)) \to H_4(\pi;\mathbb{Z})$ is the Hurewicz map.
\end{pf}
\begin{cor}\label{cor_classification_type_II}
Let $\psi \colon \pi_1 \to \pi_2$ be an epimorphism. Then the following diagram is commutative
\begin{center}
\begin{tikzpicture}
\matrix(m)[matrix of math nodes, row sep=1cm,column sep=1cm]
{
\Omega_4^{SO}(K(\pi_1,1)) & \Omega_{4}^{SO}(K(\pi_2,1)) \\
\mathbb{Z} \oplus H_4(K(\pi_1,1),\mathbb{Z}) & \mathbb{Z} \oplus H_4(K(\pi_2,1),\mathbb{Z})\\
};
\path[->,font=\scriptsize]
(m-1-1) edge node[above] {$\psi_{\ast}$} (m-1-2)
		  edge node {} (m-2-1)
(m-1-2) edge node {} (m-2-2)
(m-2-1) edge node[above] {$\operatorname{id} \oplus \psi_{\ast}$} (m-2-2);
\end{tikzpicture}
\end{center}
In this diagram $\psi_{\ast}$ occurring in $\operatorname{id} \oplus \psi_{\ast}$ is the homomorphism induced by $\psi$ on homology of $K(\pi_1,1)$ and $K(\pi_2,1)$.
\end{cor}
Corollary \ref{cor_classification_type_II} yields the following theorem.
\begin{thm}\label{classfication_odd_type}
Let $M_1$ and $M_2$ be two odd $4$-manifolds with fundamental group $\pi$. Then $Y(M_1,\psi_1)$ and $Y(M_2,\psi_2)$ are stably diffeomorphic if and only if $\sigma(M_1) = \sigma(M_2)$ and $(\psi_1)_{\ast}[M_1] = (\psi_2)_{\ast}[M_2]$.
\end{thm}

\begin{rem}
From the corollary \ref{cor_classification_type_II} it is easy to deduce that the map
$$\psi_{\ast} \colon \Omega_4^{SO}(K(F_{\alpha},1)) \to \Omega_4^{SO}(K(\pi,1))$$
is surjective if and only if $H_4(\pi;\mathbb{Z})=0$. Thus every odd $4$-manifold with fundamental group $\pi$ is stably diffeomorphic to one obtained by the surgery on loops from a manifold with free fundamental group if and only if $H_4(\pi;\mathbb{Z})=0$.
\end{rem}

Theorems \ref{classification_even_type} and \ref{classfication_odd_type} cover even and odd cases, respectively, of the proof of the main theorem \ref{main_theorem}.

\section{Thickenings of $2$-complexes}\label{section_thickenings}
Let $X$ be a finite, $2$-dimensional CW-complex. Let us denote $\pi= \pi_1(X)$. Let $\xi$ be a real vector bundle over $X$. Using the procedure described in \cite{Mazur} it is possible to construct a compact manifold $N(X,\xi)$ of dimension $2+\dim \xi$ with an inclusion $X \hookrightarrow N(X,\xi)$, which is a simple homotopy equivalence. Furthermore the map $\partial N(X,\xi) \hookrightarrow N(X,\xi)$ is $(\dim \xi -2)$-connected. The normal data of $N(X,\xi)$ is determined by the vector bundle $\xi$ in the sense that the restriction of the stable normal bundle of $N(X,\xi)$ to $X$ is stably isomorphic to $\xi$.

\begin{dfn}
Boundary of $N(X,\xi)$ is called a thickening of $X$. Thickenings will be denoted by $T(X,\xi)$, where $\xi$ denotes vector bundle used in the construction.
\end{dfn}

Thus our starting point is a $2$-dimensional CW-complex $X$ and a vector bundle $\xi$ of dimension three over $X$.

\begin{lemma}\label{lemma_splitting}
Let $\xi$ be a vector bundle of dimension at least three over a $2$-dimensional CW-complex $X$. Then there exists a splitting
$$\xi = \eta \oplus \underline{\mathbb{R}}^{\dim \xi -2},$$
where $\eta$ is s $2$-dimensional vector bundle over $X$.
\end{lemma}
\begin{pf}
The proof is an easy application of elementary obstruction theory.
\end{pf}

\begin{dfn}
Let $M$ be a compact manifold with boundary. Define the double of $M$ to be
$$D(M) = \partial (M \times I).$$
Equivalently $D(M) = M \cup_{\partial M} M$.
\end{dfn}

\begin{cor}\label{corr_double_diffeo_thickening}
Let $X$ be a $2$-dimenisional CW-complex and $\xi$ a vector bundle over $X$ of dimension three. Let $\xi = \eta \oplus \underline{\mathbb{R}}$ be a splitting as in lemma \ref{lemma_splitting}. Then there exists a diffeomorphism
$$T(X,\xi) \longrightarrow D(N(X,\eta)).$$
\end{cor}

Let's look closer at $N(X,\eta)$, where $\eta$ is a $2$-dimensional vector bundle over $X$. As was said earlier $X$ determines a presentation $\mathcal{P}$ of the fundamental group $\pi=\pi_1(X)$. While performing the construction described in \cite{Mazur} one first thickens $1$-skeleton. This produces a $\natural_n S^1 \times D^3$. Then one glues $2$-handles according to the presentation $\mathcal{P}$. Framings of the attaching circles of $2$-handles are determined by the choice of the vector bundle $\eta$, see \ref{rem_framings}. This allows us to draw a Kirby diagram of $N(X,\eta)$. In order to get a Kirby diagram of $D(N(X,\eta))$ one has to refer to \cite{gompf_stip}*{ex. 4.6.3.}. The  diagram of the double $D(N(X,\eta))$ is constructed from the one of $N(X,\eta)$ by addition of $0$-framed meridian to every attaching circle of every $2$-handle and addition of as many $3$-handles and $4$-handles as there is $1$-handles and $0$-handles in the Kirby diagram of $N(X,\eta)$. The normal $1$-type 
of $D(N(X,\eta))$ and the thickening $T(X,\eta)$ is determined by the Stiefel-Whitney class $w_2(\xi)$. 

The characteristic feature of the Kirby diagram of doubles is the fact that attaching circles of $2$-handles can be divided into two sets. One of them consists of $0$-framed meridians of the circles from the other set. This fact can be used to transform every Kirby diagram of $D(N(X,\eta))$ into a canonical form.

\begin{dfn}\label{dfn_normal_form}
Let $D$ be a Kirby diagram consisting of $n$ $1$- and $3$-handles and $2k$ $2$-handles such that the set of attaching circles of $2$-handles can be dived into two subsets $H = \{h_1,\ldots,h_k\}$ and $M=\{m_1,\ldots,m_k\}$ of equal cardinality. Every element of $m_i \in M$ is a $0$-framed meridian for $h_i \in H$ for $i=1,2,\ldots,k$. Furthermore if two circles $c_1$ and $c_2$ are geometrically linked, then there exists $1 \leq i \leq k$ such that $c_1 = m_i$ and $c_2=h_i$. Additionally we require the framing coefficient to be equal to either to $0$ or $1$. Such a Kirby diagram is said to be in the normal form.
\end{dfn}

\begin{lemma}
Kirby diagram of a double $D(N(X,\eta))$, where $\eta$ is a $2$-dimensional vector bundle over $X$, can be transformed to the normal form.
\end{lemma}
\begin{pf}
The decomposition of a set of attaching circles of $2$-handles is granted for free. To fulfill the disjointness condition use the meridians as in \cite{gompf_stip}*{prop. 5.1.4.}. Finally to adjust the framing coefficient slide elements of $H$ over corresponding meridians. All framing coefficients are even, therefore by performing appropriate number of handle slides they can be made to equal $0$.
\end{pf}

\begin{rem}\label{rem_framings}
Since $w_2(\xi)$ determines the normal $1$-type of the thickening $T(X,\xi)$, it also determines framing coefficients of its $2$-handles modulo $2$, because according to \cite{gompf_stip}*{cor. 5.7.2.} the second Stiefel-Whitney class of $T(X,\xi)$ is represented by the cocycle $c$ whose values on each $2$-handle is its framing coefficient reduced modulo $2$. If $c'$ is any other cocyle such that $c-c'$ is a coboundary, then $1$-handle twists allow to change the Kirby diagram so that the cocycle determined by the framing coefficients modulo $2$ is equal to $c'$.
\end{rem}

\section{Tietze transformations and Tietze theorem}\label{section_Tietze}

\begin{dfn}
Let $\mathcal{P}$ be a finite presentation of a group $\pi$. The following two transformations of $\mathcal{P}$ are called Tietze transformations.
\end{dfn}
\begin{enumerate}
\item Let $x$ be an arbitrary word in the alphabet generated by generators from $\mathcal{P}$ and their inverses. Transformation $T_1$ adjoins a new generator $y$ and new relator $yx^{-1}$ to $\mathcal{P}$. The transformation $T_1'$ is the reverse of $T_1$.
\item Transformations $T_2$ and $T_2'$ add or remove a relator which is a consequence of the others.
\end{enumerate}
From the definition of the Tietze transformations it follows that if two presentations can be transformed one into the other by a finite sequence of Tietze transformations, then they yield isomorphic groups. A well known theorem of Tietze assert that the opposite implication is also true.

\begin{thm}[Tietze]
Let $\mathcal{P}_1$ and $\mathcal{P}_2$ be two finite presentations of groups $\pi_1$ and $\pi_2$. Then $\pi_1$ is isomorphic to $\pi_2$ if and only if there is a finite sequence of Tietze transformations which lead from $\mathcal{P}_1$ to $\mathcal{P}_2$. 
\end{thm}

For the purpose of this paper it is be desirable to use different, but equivalent, set of transformations.
\begin{dfn}
Let $\mathcal{P}$ be a presentation of a group $\pi$. 
\begin{enumerate}
\item A trasnformation $S_1$ is a transformation of $\mathcal{P}$, which replaces any relator $r$ by $rw^{-1}sw$, where $s$ is any other relator or its inverse and $w$ is any word in the alphabet generated by generators of $\mathcal{P}$ and their inverses.

\item Transformations $S_2$ and $S_2'$ add or delete a trivial relator $e$, where $e$ denotes the identity element.
\end{enumerate}
\end{dfn}

It is easy to see that the transformations $S_1$, $S_2$ and $S_2'$ do not change the isomorphism type of the group. In fact these transformations with $T_1$ and $T_1'$ constitute a set of transformations equivalent to the set of Tietze transformations.

\section{Geometric realization of Tietze transformations}\label{section_geometric_realization}

Our aim in this paragraph is to prove theorem \ref{thm_geometric_realisation_stable_diffeo}. The proof is divided into three lemmas.

\begin{lemma}\label{transformation_T_1}
Let $\mathcal{P}_1$ and $\mathcal{P}_2$ be two presentations such that $\mathcal{P}_2$ was obtained from $\mathcal{P}_1$ by a single application of $T_1$ transformation. Let $X_1$ and $X_2$ be two $2$-dimensional CW-complexes corresponding to $\mathcal{P}_1$ and $\mathcal{P}_2$, respectively. Then there exists a diffeomorphism
$$T(X,\xi_1) \stackrel{\cong}{\longrightarrow} T(X_2,\xi_2).$$
\end{lemma}
\begin{pf}
Choose $x$ to be a word in generators of $\mathcal{P}$ and their inverses. Let $D$ be a Kirby diagram of $T(X_1,\xi_1)$. Kirby diagram of $T(X_2,\xi_2)$ is depicted on figure \ref{Kirby_diagram_T_1_start}.
\begin{figure}
\begin{center}
\includegraphics[scale=0.5]{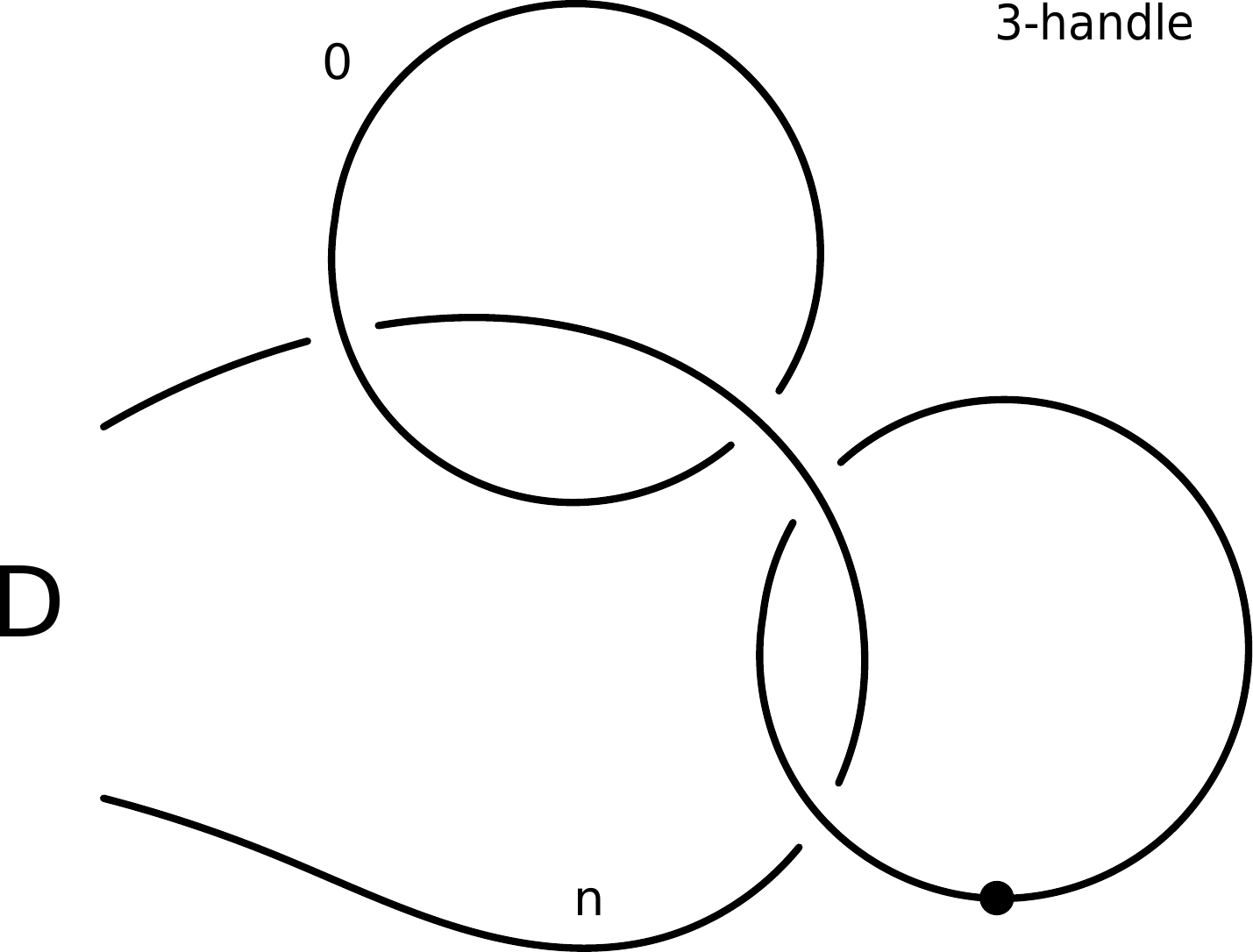}
\caption{Kirby diagram of $T(X_2,\xi_2)$.\label{Kirby_diagram_T_1_start}}
\end{center}
\end{figure}
Homotopy class of the loop, let's call it $\gamma$, which links the new $1$-handle, denoted by $\delta$, equals $x \in \pi$. Now use the $0$-framed meridian of $\gamma$ to unlink it from any $2$-handle of $D$. Slide every $1$-handle from $D$ over $\delta$ to unlink it from $D$. Repeat this process as many times as necessary to unlink $\gamma$ from any loop of $D$. After these transformations the Kirby diagram looks as the figure \ref{Kirby_diagram_T_1_step_2} shows.
\begin{figure}
\begin{center}
\includegraphics[scale=0.5]{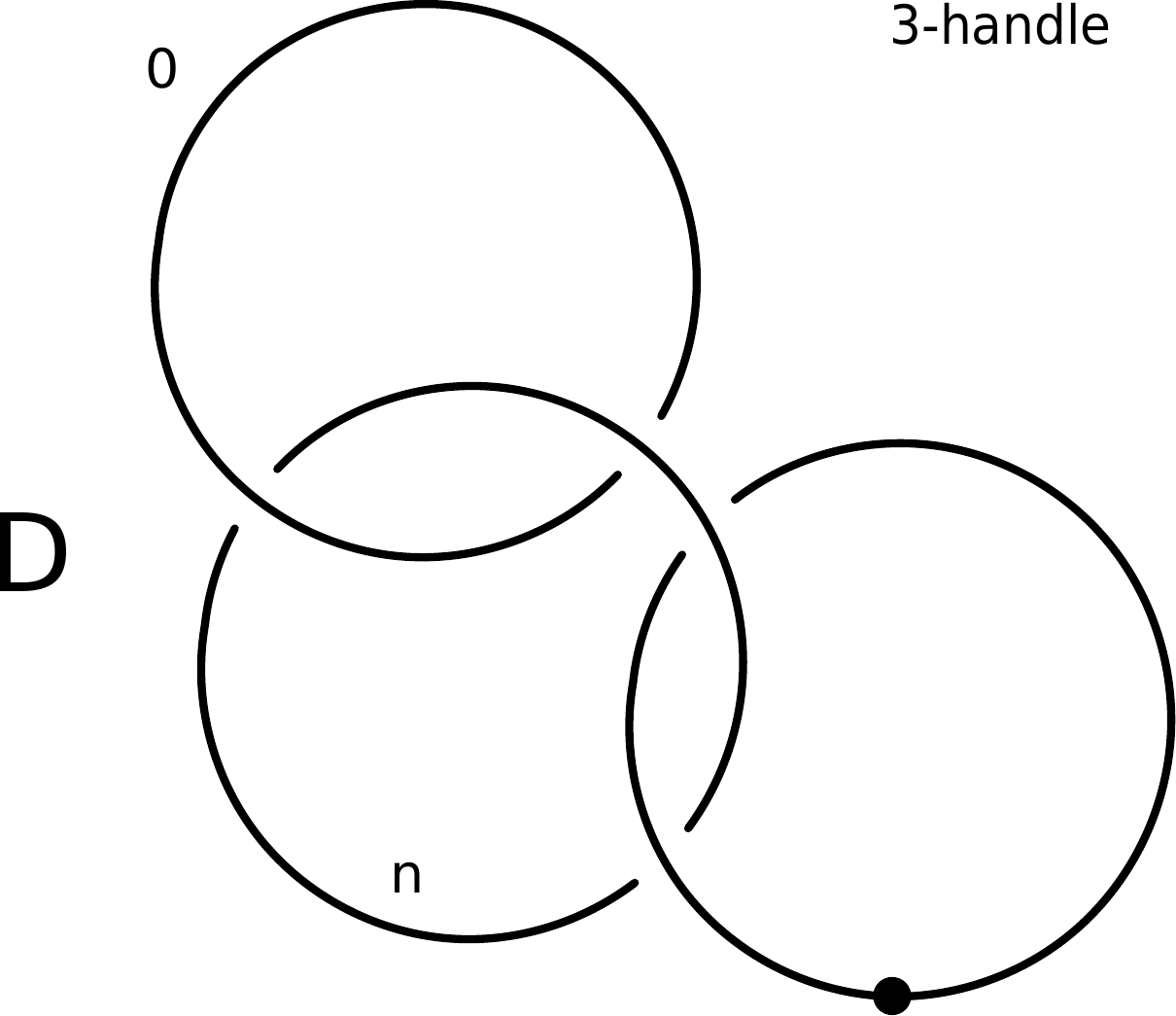}
\caption{Kirby diagram of $T(X_2,\xi_2)$ after unlinking $\gamma$ from $D$.\label{Kirby_diagram_T_1_step_2}}
\end{center}
\end{figure}
Now slide meridian of $\gamma$ over the $1$-handle and cancel it with a $3$-handle. Then cancel $\gamma$ and the $1$-handle. After performing this procedure one obtains the Kirby diagram of $T(X_1,\xi_1)$, which gives the desired conclusion.
\end{pf}

\begin{lemma}\label{transformation_S_1}
Let $\mathcal{P}_1$ and $\mathcal{P}_2$ be two presentations such that $\mathcal{P}_2$ was obtained from $\mathcal{P}_1$ by a single application of the $S_1$ transformation. Let $X_1$ and $X_2$ be two $2$-dimensional CW-complexes corresponding to $\mathcal{P}_1$ and $\mathcal{P}_2$, respectively. Then there exists a diffeomorphism
$$T(X,\xi_1) \stackrel{\cong}{\longrightarrow} T(X_2,\xi_2).$$
\end{lemma}
\begin{pf}
The part of the Kirby diagram of $T(X_1,\xi_1)$ relevant for our considerations is depicted on figure \ref{Kirby_diagram_S_1_before}.
\begin{figure}
\begin{center}
\includegraphics[scale=0.5]{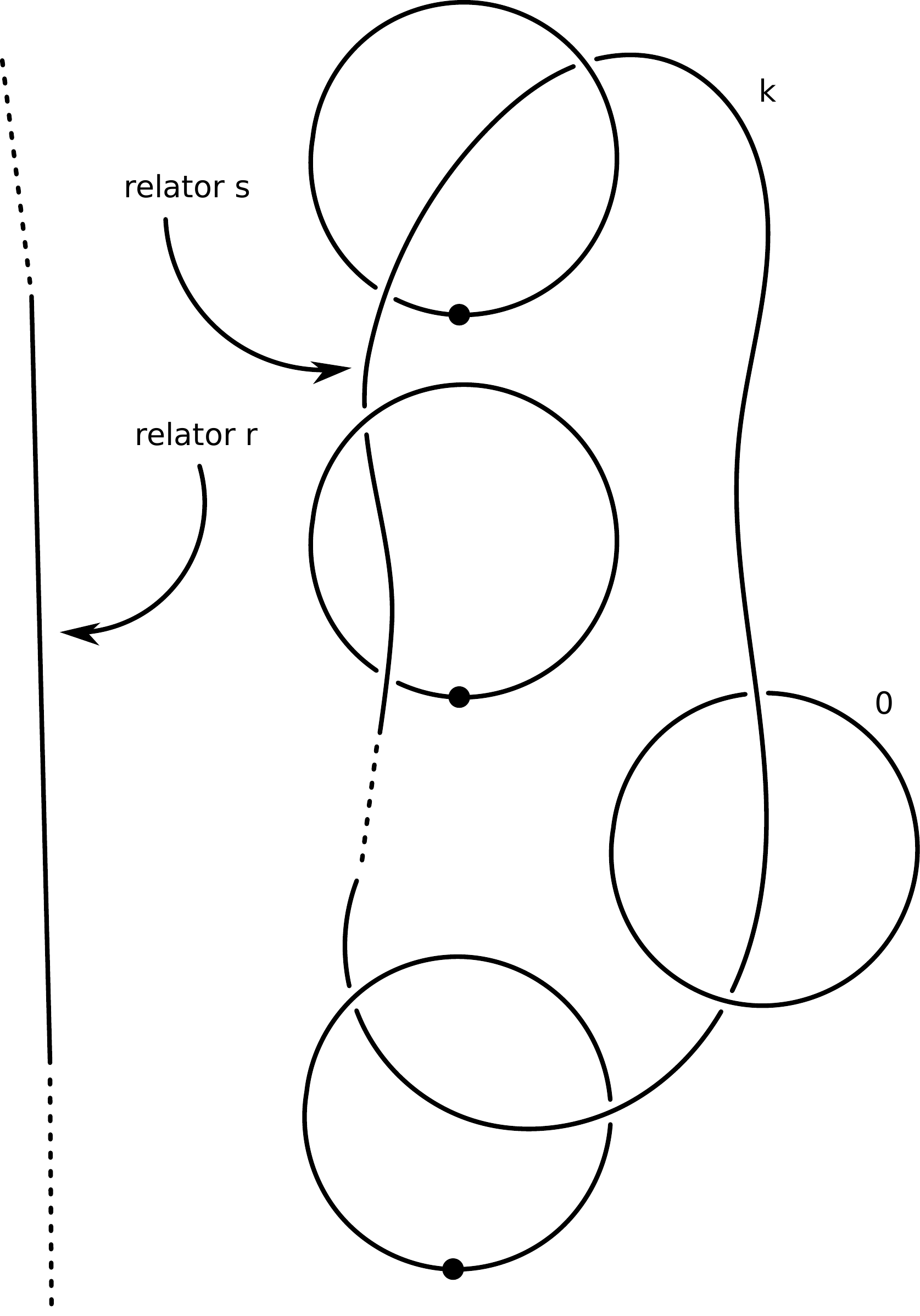}
\caption{Kirby diagram of $T(X_1,\xi_1)$.\label{Kirby_diagram_S_1_before}}
\end{center}
\end{figure}
Figure \ref{Kirby_diagram_S_1_start} shows the relevant part of Kirby diagram of $T(X_2,\xi_2)$.
\begin{figure}
\begin{center}
\includegraphics[scale=0.5]{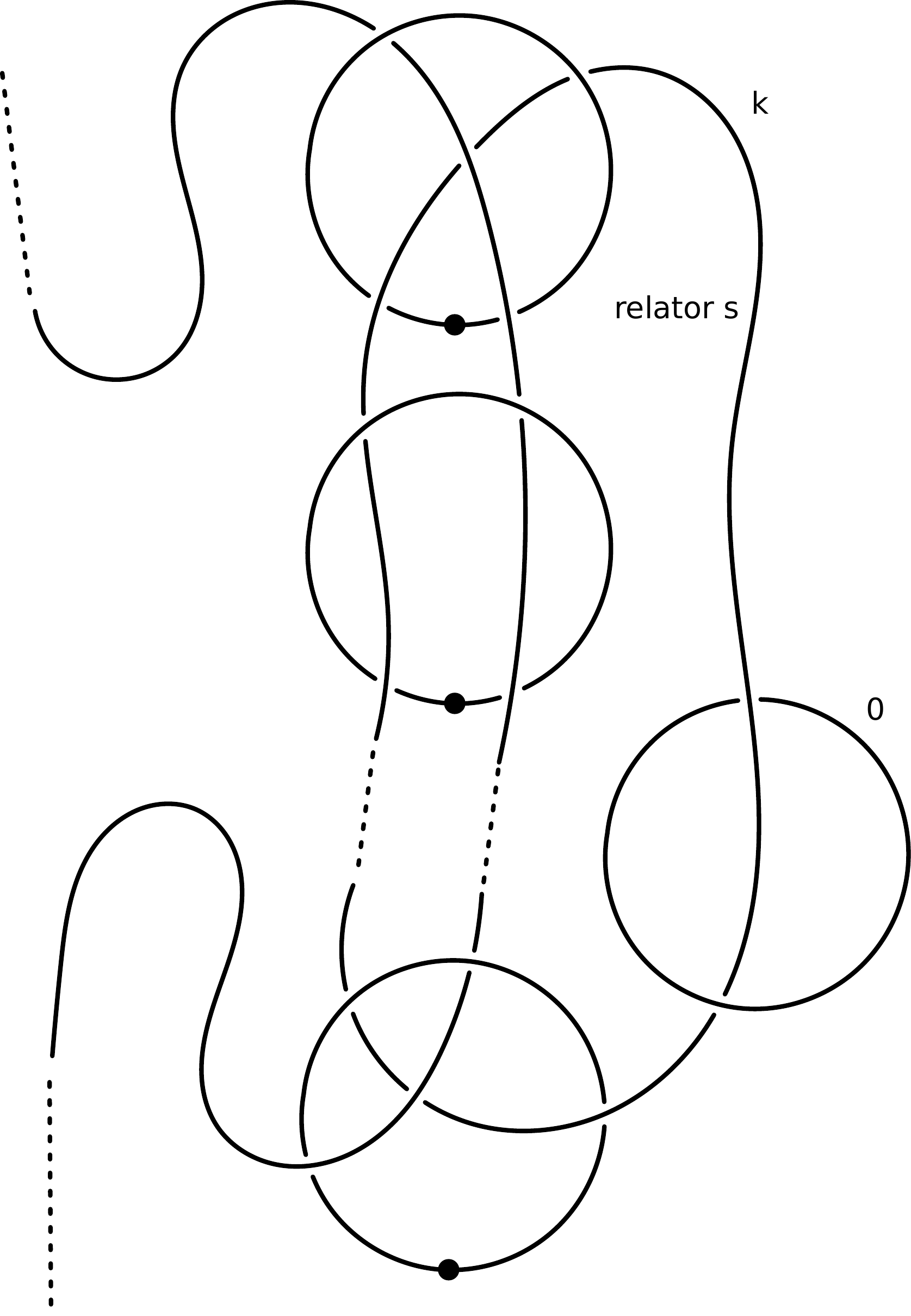}
\caption{Kirby diagram of $T(X_2,\xi_2)$.\label{Kirby_diagram_S_1_start}}
\end{center}
\end{figure}
Let $\gamma_{swrw^{-1}}$, $\mu_{swrw^{-1}}$, $\gamma_s$ and $\mu_s$ denote the attaching circles and meridians of the relators $swrw^{-1}$ and $s$, respectively.

Slide $\gamma_{swrw^{-1}}$ over $\gamma_s$ to unlink it from $1$-handles visible on the picture. Use $\mu_{s}$ to unlink $\gamma_{swrw^{-1}}$ from $\mu_s$. Now isotope $\gamma_{swrw^{-1}}$ to undo the change done by loops $w$ and $w^{-1}$. This sequence of handles slides and isotopies lead us back to the Kirby diagram of $T(X_1,\xi_1)$. 
\end{pf}

\begin{lemma}\label{transformation_S_2}
Let $\mathcal{P}_1$ and $\mathcal{P}_2$ be two presentations which differ by a single application of $S_2$ transformation. Then the exists a diffeomorphism
$$T(X_1,\xi_1) \# S^2 \times S^2 \stackrel{\cong}{\longrightarrow} T(X_2,\xi_2)$$
\end{lemma}
\begin{pf}
This is an easy consequence of the construction of $T(X_1,\xi_1)$ and $T(X_2,\xi_2)$.
\end{pf}

\begin{dfn}
The sequences of handle slides, isotopies, handle cancellations and creations as performed in the proofs of lemmas \ref{transformation_T_1}, \ref{transformation_S_1} and \ref{transformation_S_2} are called geometric Tietze transformations of type $T_1'$, $S_1'$ and $S_2'$. Geometric Tietze transformations $T_1$, $S_1$ and $S_2$ are definied as inverses of $T_1$, $S_1$ and $S_2$, respectively.
\end{dfn}

\begin{pf}[Proof ot thm. \ref{thm_main_1}]
Theorem \ref{thm_main_1} follows easily from lemmas \ref{transformation_T_1}, \ref{transformation_S_1} and \ref{transformation_S_2}.
\end{pf}

\section{The effect of a surgery on a Kirby diagram}\label{section_surgery_kirby_diagram}

Suppose that $X$ is a compact $4$-manifold with a Kirby diagram $D$. Suppose one chooses an element $x \in \pi_1(X)$ to perform a surgery on the loop $\gamma$ representing $x$. Perform the geometric $T_1$ transformation on $D$ corresponding to the following Tietze transformation
$$\left\langle x_1,\ldots,x_n | r_1,\ldots,r_m \right\rangle \longmapsto \left\langle x_1,\ldots,x_n, y | r_1,\ldots,r_m, yx^{-1} \right\rangle$$
During this transformation the following handles were added: a $1$-handle $h_1$ which links exactly once with $\gamma$, $2$-handle $h_2^1$ whose attaching circle is $\gamma$, $h_2^1$ a $0$-framed meridian for $\gamma$ and a $3$-handle $h_3$. In the new Kirby diagram $D'$ the $1$-handle $h_1$ represents a homotopy class $x$. Now perform surgery on $h_1$. This is equivalent to erasing the dot from the circle in $D'$ representing $h_1$ and giving a new $2$-handle framing equal to $0$. After the surgery the loop $\gamma$ has two $0$-framed meridians. Use one of them to unlink the other from $\gamma$ and cancel it with the $3$-handle $h_3$.
At the begging we didn't impose any condition on the framing coefficient of $\gamma$. Since after performing surgery $\gamma$ posses a $0$-framed meridian, thus only the remainder of the framing coefficient modulo 2 is relevant. The choice of even or odd framing for $\gamma$ corresponds to choosing one of the two possible framings for the loop during the surgery.

\begin{pf}[Proof of the thm. \ref{thm_diffeo_thick_surg}]
The Kirby diagram of the connected sum $\#_n S^1 \times S^3$ consists of $n$ dotted circles corresponding to $1$-handles, $n$ $3$-handles and one $4$-handle.

At the beginning of this section the description of changes made on the Kirby diagram were described. Therefore the Kirby diagram of $Y(\#_n S^1 \times S^3, \mathcal{P})$ is obtained from the Kirby diagram of $\#_n S^1 \times S^3$ by adding $2$-handles whose attaching circles represent relations of the presentation $\mathcal{P}$ and adding a $0$-framed meridian for each such circle. The only thing left to do is to perform handle slides to bring this Kirby diagram to the normal form.  

From this description it is obvious that $Y(\#_n S^1 \times S^3,\mathcal{P},w)$ is diffeomorphic to $T(X(\mathcal{P}),\xi)$, where $X(\mathcal{P})$ denotes the $2$-dimensional CW complex corresponding to the presentation $\mathcal{P}$ and $\xi$ denotes the appropriate bundle.
\end{pf}

\noindent Department of Mathematics and Computer Science, Adam Mickiewicz University in Poznan, ul. Umultowska 87, 61-614 Poznan, Poland,

\noindent E-mail address: \textit{politarw@amu.edu.pl.}

\end{document}